\documentclass[a4paper]{article}
\usepackage{amsmath}
\usepackage{amssymb}
\newcommand{\be}{\begin{equation}}
\newcommand{\ee}{\end{equation}}
\newcommand{\ra}{\rightarrow}

\newcommand{\bea}{\begin{eqnarray}}
\newcommand{\eea}{\end{eqnarray}}
\newcommand{\beas}{\begin{eqnarray*}}
\newcommand{\eeas}{\end{eqnarray*}}

\newcommand{\R}{\mathbb{R}}
\newcommand{\C}{\mathbb{C}}
\newcommand{\F}{\mathbb{F}}
\newcommand{\1}{{\mathbf{1}}}

\newcommand{\A}{{\cal A}}

\def\Ci{C^\infty}
\def\M{\mathfrak{M}}
\def\S{\mathfrak{S}}
\def\di{\mathfrak{D}}
\def\m{\mathfrak{m}}

\def\c{\mathfrak{c}}

\mathchardef\za="710B  
\mathchardef\zb="710C  
\mathchardef\zg="710D  
\mathchardef\zd="710E  
\mathchardef\zve="710F 
\mathchardef\zz="7110  
\mathchardef\zh="7111  
\mathchardef\zvy="7112 
\mathchardef\zi="7113  
\mathchardef\zk="7114  
\mathchardef\zl="7115  
\mathchardef\zm="7116  
\mathchardef\zn="7117  
\mathchardef\zx="7118  
\mathchardef\zp="7119  
\mathchardef\zr="711A  
\mathchardef\zs="711B  
\mathchardef\zt="711C  
\mathchardef\zu="711D  
\mathchardef\zvf="711E 
\mathchardef\zq="711F  
\mathchardef\zc="7120  
\mathchardef\zw="7121  
\mathchardef\ze="7122  
\mathchardef\zy="7123  
\mathchardef\zf="7124  
\mathchardef\zvr="7125 
\mathchardef\zvs="7126 
\mathchardef\zf="7127  
\mathchardef\zG="7000  
\mathchardef\zD="7001  
\mathchardef\zY="7002  
\mathchardef\zL="7003  
\mathchardef\zX="7004  
\mathchardef\zP="7005  
\mathchardef\zS="7006  
\mathchardef\zU="7007  
\mathchardef\zF="7008  
\mathchardef\zW="700A  

\newcommand{\ep}{\hfill$\Box$}
\newcommand{\bp}{\textit{Proof.-} }

\begin{document}

\title{Isomorphisms of algebras of smooth functions\\ revisited
\footnote{This work was supported by KBN, grant No 2 P03A 020 24.}}
\author{J. Grabowski}\maketitle

\newtheorem{re}{Remark}
\newtheorem{theo}{Theorem}
\newtheorem{prop}{Proposition}
\newtheorem{lem}{Lemma}
\newtheorem{cor}{Corollary}
\newtheorem{ex}{Example}

\begin{abstract}
It is proved that isomorphisms between algebras of smooth
functions on Hausdorff smooth manifolds are implemented by
diffeomorphisms. It is not required that manifolds are connected
nor second countable nor paracompact. This solves a problem stated
by A.~Weinstein. Some related results are discussed as well.



\end{abstract}

\section{Introduction}
Choose $\F$ to be either $\R$ or $\C$. The following looks
familiar to most mathematicians.

\begin{theo}\label{t1} Every algebra isomorphism $\zF:\Ci(M_1;\F)
\ra\Ci(M_2;\F)$ between the associative algebras of all $\F$-valued smooth
functions on Hausdorff smooth finite-dimensional manifolds $M_1$ and $M_2$
is the pullback by a smooth diffeomorphism $\zf:M_2\ra M_1$.
\end{theo}

However, as it was pointed out to us by A.~Weinstein, the standard
proofs of this fact available in the literature strongly use the
additional requirement that the manifolds are second countable.
This is because they use the interpretation of points of such
manifolds as multiplicative functionals on the corresponding
algebras of functions (this result is sometimes called "Milnor's
exercise", cf. \cite[p.~11]{MS}), that has been proved for second
countable manifolds (see e.g. \cite[Prop.~3.5.]{JG} or
\cite[Suppl. 4.2C]{AMR}). Even if we assume that the manifolds are
paracompact such proofs work only when the number of the connected
components of the manifolds (e.g. discrete sets) is not bigger
than the cardinality of the reals.

\medskip
Of course, a similar problem occurs when we deal with algebras
$C(X;\F)$ of all continuous instead of smooth functions. The fact
that, for $X_i$ being compact (respectively, completely regular
and first countable) topological spaces, $i=1,2$, the associative
algebras $C(X_1;\R)$ and $C(X_2;\R)$ are isomorphic if and only if
$X_1$ and $X_2$ are homeomorphic was proved already in 1937 by
Gel'fand and Kolmogoroff \cite{GK} (see also \cite{St}). In the
compact case the authors used also the characterization of points
of these spaces as multiplicative functionals (or, equivalently,
one-codimensional ideals) of the corresponding algebras of
continuous functions. In the general case, an identification of
the space of all maximal ideals with the Stone-\v Cech
compactification $\zb X$, together with some properties of $\zb X$
for first countable $X$, was used. Note that the theorems in
\cite{GK} have an existential character and the form of the
algebra isomorphism is not given.

\medskip
The class of completely regular topological spaces $X$ such that
every one-codimensional ideal in $C(X;\R)$ is of the form $p^*=\{
f\in C(X;\R): f(p)=0\}$ for a certain $x\in X$ appeared in
\cite{He1} under the name \textit{Q-spaces} (there are various
equivalent definitions). Now the name \textit{realcompact spaces}
is commonly used. Let us call ideals of the form $p^*$
\textit{fixed} and the other one-codimensional ideals
\textit{free}. In this language, the space $X$ is realcompact if
all one-codimensional ideals in $C(X;\R)$ are fixed. Smooth
manifolds $M$ with the analogous property of the algebra
$\Ci(M;\R)$ are called \textit{smoothly realcompact}. For a survey
of results on properties of realcompact and smoothly realcompact
spaces we refer to \cite[Ch.~3.11]{En}, \cite[Ch.~8]{GJ}, and
\cite[Ch.~IV]{KM}.

Note that the problem of realcompactness of discrete sets reduces
to the problem of $\zs$\textit{-measurability} of their
cardinalities (see Definition 1 and the remarks thereafter). In
general, paracompact spaces are realcompact if and only if the
cardinalities of all their closed discrete subsets are not
$\zs$-measurable \cite{Ka,He2} (see also \cite[5.5.10]{En}). Since
connected paracompact Hausdorff smooth manifolds are second
countable, the maximal cardinality of closed discrete subsets of a
paracompact Hausdorff manifold $M$ which is not second countable
equals the cardinality of the set of all connected components of
$M$. We will prove the following.

\begin{theo}\label{t2} A paracompact smooth manifold $M$ is realcompact if and
only if the cardinality $\m$ of the set of all components of $M$
is not $\zs$-measurable.
\end{theo}

This means that one cannot identify points with one-codimensional
ideals in $\Ci(M;\R)$ for a paracompact Hausdorff smooth manifold
$M$ with $\zs$-measurable cardinality of components, so one cannot
apply directly the standard proof of the form of isomorphisms of
the algebras of smooth functions for such manifolds. Of course,
one can try to adapt the proof of Gel'fand and Kolmogoroff
\cite{GK}, but there are several delicate points there.

\medskip
There are very few papers on differentiable manifolds which are
not assumed to be paracompact. This is because main tools like the
partition of unity are not available in that case. Our aim in this
note is to prove Theorem \ref{t1} in full generality. The trick is
that, to characterize points, we use not all one-codimensional
ideals but a natural subclass of them. We present also a few
related results. In particular, we get a short proof of the
Gel'fand-Kolmogoroff result \cite{GK} complemented by a
description of the form of isomorphisms and without use of the
Stone-\v Cech compactification.

\medskip
After publishing our proof in the arxiv, we found the preprint
\cite{Mrc}, where Theorem \ref{t1} has been proved (first for
paracompact manifolds, then, in a new version of the preprint, in
general) by different methods using \textit{characteristic
sequences of functions} instead of the characterization of
multiplicative functionals on the algebras of smooth functions.

\medskip
All smooth manifolds in this note are assumed to be Hausdorff and
finite-dimensional if not otherwise stated.

\section{Smoothly realcompact manifolds}

For the convenience of the reader let us start with recalling some notions
from Set Theory.

\medskip\noindent
\textbf{Definition 1.} By a $\{ 0,1\}$-\textit{valued}
$\zs$\textit{-measure} on a set $X$ we mean a countably additive
function $\zm$ defined on the family of all subsets of $X$, and
assuming only the values $0$ and $1$. We call such a measure
\textit{free} if $\zm(\{ x\})=0$ for all $x\in X$ and
\textit{trivial} if $\zm\equiv 0$. A cardinal $\m$ we call
$\zs$\textit{-measurable} if a set $X$ of cardinality $\m$ admits
a $\{ 0,1\}$-valued $\zs$-measure which is free and nontrivial.

\bigskip\noindent
\textbf{Remarks.} \begin{enumerate} \item Sometimes in the
literature $\zs$-measurable cardinals defined above are called
just measurable. We decided to distinguish $\zs$-measurability
from the notion of \textit{measurabilty} of cardinals which is
used nowadays in Set Theory: an uncountable cardinal $\m$ is
\textit{measurable} if there exists an $\m$-complete nonprincipal
(free) ultrafilter over $\m$ or, equivalently, if there exists a
non-trivial $\{ 0,1\}$-valued measure on $\m$ which is
$\zk$-additive for all $\zk<\m$ (cf. \cite[Ch.~5]{Je}). Since the
smallest $\zs$-measurable cardinal is measurable (\cite[Lemma
27.1]{Je}) the problem of existence of measurable cardinals is
equivalent to the problem of existence of $\zs$-measurable
cardinals (and far from being solved).

\item It is obvious that $\aleph_0$ is not $\zs$-measurable, so,
according to a theorem by S.~Ulam \cite{Ul}, $\c=2^{\aleph_0}$ is
not $\zs$-measurable. Since each $\{ 0,1\}$-valued $\zs$-measure
is $\m$-additive for every non-$\zs$-measurable $\m$ (cf.
\cite[12.3]{GJ}), each $\{ 0,1\}$-valued $\zs$-measure $\zm$ is
$\c$-additive, i.e. $\zm(\bigcup_{\zg<\c}S_\zg)=1$ implies that
$\zm(S_\zg)=1$ for certain $\zg<\c$. Here and further we identify
cardinal numbers with the smallest ordinals with the same
cardinality.
\end{enumerate}

\bigskip
We will prove a result slightly more general than Theorem \ref{t2}
(algebras of complex-valued functions are included).

\begin{theo}\label{tn} Let $M$ be a paracompact smooth manifold $M$. Then
there is a free one-codimensional ideal in $\Ci(M;\F)$ if and only
if the cardinality $\m$ of the set of all connected components of
$M$ is $\zs$-measurable.
\end{theo}
\bp Let $C=\{ C_\zg:\zg<\m\}$ be the set of all connected components of
$M$.

($\Leftarrow$) Suppose $\m$ is $\zs$-measurable and let $\zm$ be a
$\{ 0,1\}$-valued free and nontrivial $\zs$-measure on
$X=\{\zg<\m\}$. By the above remark, $\zm$ is $\c$-additive.
Consider a choice $\{ p_\zg\in C_\zg:\zg<\m\}$ of points of the
components of $M$ and put
$$J=\{ f\in\Ci(M):\zm(\{\zg<\m: f(p_\zg)=0\})=1\}.
$$
It is easy to see that $J$ is a proper ideal in $\A=\Ci(M)$.
Moreover, $J$ is one-codimensional in $\A$. Indeed, for any
$g\in\A$ consider the partition of $X$ consisting of subsets
$V_r=g^{-1}(\{ r\})$, $r\in\R$. Since $\zm$ is $\c$-additive,
$\zm(V_{r_0})=1$ for a certain $r_0\in\R$. But this means that
$(g-r_0\cdot \1_M)\in J$, so $J$ is one-codimensional. Finally,
that $J$ is free follows from the fact that $\zm$ is free.

($\Rightarrow$) Suppose that there is a free one-codimensional
ideal $J$ in $\A=\Ci(M)$. We define a $\{ 0,1\}$-valued nontrivial
function $\zm$ defined on the subsets of $X=\{\zg<\m\}$ by
$$\zm(S)=1\ \Leftrightarrow\ \exists f_S\in J\ [S=\{\zg<\m:f_S^{-1}(0)
\cap C_\zg\ne\emptyset\}].
$$
We will show that $\zm$ is countably additive. Observe first that
$\zm(S)=\zm(S')=1\ \Rightarrow S\cap S'\ne\emptyset$. Indeed, if
one had $S\cap S'=\emptyset$ then the zero-sets of $f_S$ and
$f_{S'}$ are disjoint, so $\vert f_S\vert^2+\vert f_{S'}\vert^2$
would be a nowhere-vanishing, so invertible, function in $J$. We
have to show that for every partition $X=\bigcup_{n=1}^\infty X_n$
by pairwise disjoint subsets there is $n_0$ such that
$\zm(X_{n_0})=1$. For, take $f\in\A$ assuming only natural values
such that $f_{\mid C_\zg}\equiv n$ for $\zg\in X_n$. Since $J$ is
one-codimensional, $f-r\cdot\1_M\in J$ for a certain $r\in\R$. It
is clear that $r$ has to be natural, say $r=n_0$, so
$\zm(X_{n_0})=1$ by definition.

Finally, we will show that the $\zs$-measure $\zm$ is free. In the
other case we would have $\zm(\{\zg_0\})=1$ for a certain
$\zg_0<\m$. This means that every function from $J$ has zeros in
$C_{\zg_0}$, i.e. $J_{\zg_0}=J\cap\Ci(C_{\zg_0};\F)$ is a
nontrivial, thus one-codimensional, ideal in $\Ci(C_{\zg_0};\F)$.
Here, of course, we understand $\Ci(C_{\zg_0};\F)$ as the
subalgebra in $\Ci(M;\F)$ consisting of functions with support in
$C_{\zg_0}$. Since $C_{\zg_0}$ is paracompact and second
countable, we are in the standard case and $J_{\zg_0}=p^*$ for
some $p\in C_{\zg_0}$. Consequently, $J=p^*$; a contradiction. \ep

\section{Distinguished ideals and isomorphisms: smooth case}

Let $\A$ be an associative commutative algebra with unit $\1$ over
a field $\mathfrak{k}$ and let $\M(\A)$ be the set of all
one-codimensional ideals in $\A$ (or, equivalently, of all
multiplicative functionals $m:\A\ra{\mathfrak{k}}$).

\medskip\noindent
\textbf{Definition 2.} An ideal $I\in\M(\A)$ is called
\textit{distinguished} if
$$I\nsubseteq\bigcup_{J\in\M(\A),J\ne I}J,
$$
i.e., if there is $f\in I$ which belongs to no other one-codimensional
ideal of $\A$.

\medskip
Denote the set of all distinguished ideals of $\A$ by $\di(\A)$. On
$\di(\A)$ we introduce the topology (sometimes called the \textit{Stone
topology}) by defining the closure $cl(S)$ of $S\subset\di(\A)$ as
consisting of those $I\in\di(\A)$ which include $\bigcap_{J\in S}J$ (cf.
\cite{GK,St}).
\begin{theo}\label{t3} Let $M$ be a Hausdorff finite-dimensional smooth manifold
and let $\A=\Ci(M;\F)$. Then,
\be\label{co} M\ni p\mapsto p^*=\{
f\in\A:f(p)=0\}\in\M(\A)
\ee
establishes a homeomorphism of $M$ onto $\di(\A)$.
\end{theo}
\bp We will show first that (\ref{co}) establishes a one-to-one
correspondence between $M$ and $\di(\A)$.

Let $I$ be a distinguished ideal of $\A$ and suppose that $f\in\A$
does not belong to any other one-codimensional ideal. Since $f$
belongs to a single one-codimensional ideal and the ideals $p^*$,
for $p\in M$, are clearly one-codimensional and pairwise
different, $f$ vanishes at not more than one point. But $f$ has to
vanish at a point, say $p$, since otherwise it is invertible in
$\A$. This implies that $I\subset p^*$, so $I=p^*$.

Conversely, it is easy to see that the ideals $p^*$, $p\in M$, are
distinguished. For, it suffices to consider a smooth function $f$
vanishing exactly at $p$ (e.g. to take locally, in a coordinate
chart $(U,x)$ centered at $p$, the function $\sum_i x_i^2$ and to
extend it smoothly to a positive function outside $U$). If $f$
were a member of any other one-codimensional ideal, say $J$, then
we would have some $g\in J$ not vanishing at $p$ and the function
$\vert f\vert^2+\vert g\vert^2$ would be an invertible member of
$J$.

Having established the identification of $M$ with $\di(\A)$ we
will finish with showing, completely analogously to \cite{GK},
that it identifies also the topologies on $M$ and $\di(\A)$.

Indeed if $S\subset M$ and $p\in cl_M(S)$ then, due to the
continuity of smooth functions, any function vanishing on $S$ has
to vanish at $p$. Conversely, if $p\notin cl_M(S)$ then we can
find a coordinate neighbourhood $U$ of $p$ not intersecting $S$
and a bump function $f\in\A$ with support in $U$ and $f(p)=1$.
Then $f\in\bigcap_{q\in S}q^*$ but $f\notin p^*$. \ep

\medskip
It is obvious that the property "to be a distinguished ideal" is a
purely algebraic property respected by algebra isomorphisms, so we
can proceed as in the standard case.

\medskip\noindent
\textit{Proof of Theorem 1.- } Denote $\A_i=\Ci(M_i;\F)$, $i=1,2$. The
algebra isomorphism $\zF:\A_1 \ra\A_2$ induces a bijection $\di(\A_2)\ni
I\mapsto\zF^{-1}(I)\in\di(\A_1)$ and, in view of Theorem \ref{t3}, a
bijection $\zf:M_2\ra M_1$ such that
\be\label{bi}
\zF(f)(p)=0\ \Leftrightarrow\ f(\zf(p))=0
\ee
for all $f\in\A_1$, $p\in M_2$. Now, since $f-f(\zf(p))\cdot\1_{M_1}$
vanishes at $\zf(p)$,
$$\zF(f-f(\zf(p))\cdot\1_{M_1})=\zF(f)-f(\zf(p))\cdot\1_{M_2}
$$
vanishes at $p$ according to (\ref{bi}), so $\zF(f)(p)=f(\zf(p))$ for all
$f\in\A_1$, $p\in M_2$, i.e. $\zF$ is the pullback by $\zf$. It remains to
prove that $\zf$ is a diffeomorphism. To be able to check smoothness in
local charts, let us show first that it is a homeomorphism. According to
Theorem \ref{t3},
$$ p\in cl_{M_2}(S)\Leftrightarrow \bigcap_{q\in
S}q^*\subset p^*\Leftrightarrow\bigcap_{q\in S}\zf(q)^*\subset \zf(p)^*\
\Leftrightarrow \zf(p)\in cl_{M_1}(\zf(S))
$$
that proves the continuity. Since $\zf^{-1}$ is continuous as well, we
conclude that $\zf$ is a homeomorphism. Now we can use the fact that in a
neighbourhood of any point $M_2$ we can use certain compactly supported
$x^1,\dots,x^n\in\A$ as local coordinates. The functions
$x^1\circ\zf,\dots,x^n\circ\zf$ are smooth on $M_1$, so $\zf$ is smooth.
Similarly, $\zf^{-1}$ is smooth, so $\zf$ is a diffeomorphism. \ep

\medskip\noindent
\textbf{Remark.\ } In Theorem 1, isomorphisms cannot be replaced
by homomorphisms, even surjective ones. A simple example is as
follows. Take $M$ which admits a free one-codimensional ideal $J$
of $\A=\Ci(M;\R)$ and the canonical projection $\zF:\A\ra\R=\A/J$.
We can consider $\R$ as the algebra of smooth functions on a
single point but $\zF$ cannot be the pullback of an embedding of
this point into $M$, since $J$ is free.

\section{Isomorphisms: general case}

It is completely obvious that a main part of the above proof
remains valid if we replace the algebras of smooth functions with
certain unitary subalgebras $\S_i$ of the algebras $C(X_i,\F)$ of
all $\F$-valued continuous functions on topological spaces $X_i$
such that
\begin{enumerate}
\item if $f\in\S_i$ is nowhere vanishing, then $f^{-1}\in\S_i$,
\item for every $p\in X_i$ and every open neighbourhood $U$ of $p$
there is $g\in\S_i$, $g:X_i\ra[0,1]$, $supp(g)\subset U$, and such that
$g(p')=1$ if and only if $p'=p$,
\end{enumerate}
$i=1,2$.
Note that if $g\in\S_i$ is as above, $1-g$ vanishes exactly at
$p$, so the above conditions ensure $\S_i$-regularity of $X_i$ and the
fact that $p^*$ are distinguished ideals in $\S_i$, $i=1,2$. Algebras of
continuous functions satisfying the above conditions we will call {\it
distinguishing}. Thus we get the following.

\begin{theo}\label{t5} Let $\S_i$ be a distinguishing algebra of $\F$-valued
continuous functions on a topological space $X_i$, $i=1,2$. Then, every
algebra isomorphism $\zF:\S_1 \ra \S_2$ is the pullback by a homeomorphism
$\zf:X_2\ra X_1$.
\end{theo}
\begin{cor}\label{t4} \cite{GK} If $X_i$, $i=1,2$, are first countable
completely regular topological spaces then every algebra isomorphism
$\zF:C(X_1;\F) \ra C(X_2;\F)$ is the pullback by a homeomorphism
$\zf:X_2\ra X_1$.
\end{cor}
\bp It suffices to prove that the algebra $C(X;\F)$ is distinguishing for
any first countable completely regular $X$. For $p\in X$ and an open
neighbourhood $U$ of $p$, we construct a continuous function
$g:X\ra[0,1]$, $supp(g)\subset U$, and such that $g(p')=1$ if and only if
$p'=p$ as follows. Take a countable basis $\{ U_n: n=1,2,\dots\}$ of the
topology at $p$, consisting of open sets contained in $U$, and, using the
complete regularity, take functions $g_n\in C(X;[0,1])$ such that
$supp(g_n)\subset U_n$ and $g_n(p)=1$. Then
$$g=\sum_{n=1}^\infty\frac{1}{2^n}g_n
$$
is the required function. \ep

\medskip\noindent
{\bf Remark.} The above results easily imply that Theorem \ref{t1}
remains valid for manifolds of class $C^k$, the algebras of
functions of class $C^k$, and diffeomorphisms of class $C^k$,
$k=0,1,\dots,\infty$. It is also valid for infinite-dimensional
manifolds of various types (e.g. modelled on Banach spaces or just
convenient vector spaces \cite[Ch.~VI]{KM}) if only the existence
of appropriate bump functions is ensured, i.e. if there are smooth
functions $f$ with supports in a given neighbourhood of $0$ in the
model topological vector space $X$ and such that $f(x)=1$ if and
only if $x=0$. This is true, for instance, for Banach spaces
admitting an equivalent smooth norm (e.g. Hilbert spaces or
$L^p(\R)$ for $p$ even). We will not discuss these problems here.
For the questions of existence of smooth bump functions and
partitions of unity on infinite-dimensional manifolds we refer to
\cite{DGZ}, \cite[Suppl.~5.5]{AMR}, and \cite[Ch.~III]{KM}.

\section{Acknowledgement.}

The author is grateful to Z.~Adamowicz, S.~Kwapie\'n, and
P.~W.~Michor for providing helpful comments and bibliographical
advice.

\noindent Janusz GRABOWSKI\\Polish Academy of Sciences\\Institute of
Mathematics\\\'Sniadeckich 8\\P.O. Box 21\\00-956 Warsaw,
Poland\\Email: jagrab@impan.gov.pl\\\\

\begin{thebibliography}{Dillo 83}
\bibitem[AMR88]{AMR} Abraham, R.; Marsden, J.~E.; Ratiu, T.:
\textit{Manifolds, Tensor Analysis, and Applications}, Springer-Verlag,
New York 1988.
\bibitem[DGZ93]{DGZ} Deville, R.; Godefroy, G.; Zizler, V.~E.:
\textit{Smoothness and Renorming in Banach Spaces}, Pitman Monographs and
Surveys in Pure and Applied Mathematics 64, Longman, John Wiley, London,
New York, 1993.
\bibitem[Eng89]{En} Engelking, R.: \textit{General Topology, revised and
completed edition}, Heldermann-Verlag, Berlin 1989.
\bibitem[GeK37]{GK} Gel'fand, I.; Kolmogoroff, A.: {On rings of continuous
functions on topological spaces}, C. R. (Dokl.) Acad. Sci. URSS,
\textbf{22} (1939), pp. 11-15.
\bibitem[GiJ60]{GJ} Gillman, L.; Jerison, M.: \textit{Rings of Continuous
Functions}, New York 1960.
\bibitem[Gra78]{JG} Grabowski, J.: \textit{Isomorphisms and ideals of the Lie
algebras of vector fields}, Invent. math., \textbf{50} (1978), pp. 13-33.
\bibitem[Hew48]{He1} Hewitt, E.: \textit{Rings of real valued continuous
functions I}, Trans. Amer. Math. Soc. \textbf{64} (1948), pp. 45-99.
\bibitem[Hew50]{He2} Hewitt, E.: \textit{Linear functionals on spaces of
continuous functions}, Fund. Math. \textbf{37} (1950), pp. 161-189.
\bibitem[Je78]{Je} Jech, T.: \textit{Set Theory}, Academic Press, New York
1978.
\bibitem[Kat51]{Ka} Kat\v etov, M.: \textit{Measures in fully normal spaces},
Fund. Math. \textbf{38} (1951), pp. 73-84.
\bibitem[KM97b]{KM} Kriegl, A.; Michor, P.~W.: \textit{The Convenient Setting
of Global Analysis}, Math. Surv. and Monogr. \textbf{53}, Amer. Math.
Soc., Providence 1997.
\bibitem[MiS74]{MS} Milnor, John W.; Stasheff, James D.:
\textit{Characteristic Classes}, Annals of Mathematics Studies, No. 76,
Princeton University Press, Princeton, N. J.; University of Tokyo Press,
Tokyo, 1974.
\bibitem[Mr\v c03]{Mrc} Mr\v cun, J.: \textit{On isomorphisms of algebras of
smooth functions}, arXiv: math.DG/0309179.
\bibitem[Sto37]{St} Stone, M.~H.: \textit{Applications of the theory of
Boolean rings to general topology}, Trans. Amer. Math. Soc. \textbf{41}
(1937), pp. 375-481.
\bibitem[Ula30]{Ul} Ulam, S.: \textit{Zur Masstheorie in der allgemainen
Mengenlehre}, Fund. Math. \textbf{16} (1930), pp. 140-150.

\end{thebibliography}
\end{document}